\documentclass[a4, 12pt]{amsart}
\usepackage[mathscr]{eucal}
\usepackage{amssymb}
\usepackage{latexsym}
\usepackage{amsthm}
\theoremstyle{plain}

\makeatletter

\title[Estimates for  lower order eigenvalues ]{ Estimates for  lower order eigenvalues  \\
of a clamped plate problem*
}
\author{Qing-Ming Cheng, Guangyue Huang and\ Guoxin Wei} 

\begin{document}
\maketitle

\begin{abstract}
\noindent For a bounded domain $\Omega$ in a complete Riemannian
manifold $M^n$, we study estimates for  lower order  eigenvalues of
a clamped plate problem. We obtain  universal inequalities for lower 
order eigenvalues. We would like to remark  that our results are sharp.
\end{abstract}



\footnotetext{* Research partially Supported by a Grant-in-Aid for
Scientific Research from JSPS.}

\section *{1. Introduction}

Let $\Omega$ be a bounded domain in an $n$-dimensional complete Riemannian manifold $M^n$.
 Assume that $\Gamma_i$ is the $i^{th}$
eigenvalue of  {\it a clamped plate problem}, which describes characteristic vibrations
of a clamped plate:

\begin{equation*}
  {\begin{cases}
     \Delta^2 u = \Gamma  u,& \ \ {\rm in} \ \ \Omega ,\\
     u=\dfrac{\partial u}{\partial \nu}=0 , & \ \ {\rm on}  \ \ \partial \Omega,
     \end{cases}}
     \eqno{(1.1)}
\end{equation*}
where $\Delta$ is the Laplacian   on $M^n$ and
$\nu$ denotes the outward unit normal to the boundary $\partial \Omega$. It is well known that this
problem has a real and discrete spectrum

$$0<\Gamma_1\leq\Gamma_2\leq\cdots\leq\Gamma_k\leq\cdots\nearrow
+\infty,$$ where each $\Gamma_i$ has finite multiplicity which is
repeated according to its multiplicity.

In this paper, we  do not introduce  results on universal inequalities for higher order eigenvalues
of the clamped plate problem. The readers who are interested
in it can see the references \cite{CQ},  \cite{[CIM1]}, \cite{[CIM2]}, \cite{H}, \cite{[PPW]} and \cite{[WX]}.
Since for lower order eigenvalues of the clamped plate problem, one can obtain
better universal inequalities for eigenvalues. We will focus our mind on the investigation of 
lower order eigenvalues of the clamped plate problem.

When $M^n$  is an $n$-dimensional Euclidean space, 
for lower order eigenvalues of the clamped plate problem (1.1),
Ashbaugh \cite{[A]}  announced  the following two universal inequalities
without proofs.  Cheng, Ichikawa and Mametsuka \cite{[CIM1]} have given
their proofs.
$$
\sum_{i=1}^n(\Gamma_{i+1}^{\frac{1}{2}}-\Gamma_1^{\frac{1}{2}})\leq
4\Gamma_1^{\frac{1}{2}},\eqno{(1.2)}
$$

$$
\sum_{i=1}^n(\Gamma_{i+1}-\Gamma_1)\leq 24\Gamma_1.\eqno{(1.3)}
$$

When $M^n$ is a  general complete Riemannian manifold
other than the Euclidean space,  it is natural to consider
the following problem:

\par\bigskip \noindent{\bf Problem}. {\it Let $M^n$ be an
$n$-dimensional complete Riemannian manifold and $\Omega$  a
bounded domain in $M^n$.  Whether can one  obtain a universal inequality  for
lower order eigenvalues, which are analogous to {\rm (1.2)},   of the clamped plate problem?}\\

In this paper, we will answer  the problem and prove the following results:

\par\bigskip \noindent{\bf Theorem 1}. {\it Let $\Omega$ be a
bounded domain in an $n$-dimensional complete Riemannian manifold
$M^n$. For the lower order eigenvalues of the clamped plate problem:
\begin{equation*}
  {\begin{cases}
     \Delta^2 u = \Gamma  u,& \ \ {\rm in} \ \ \Omega ,\\
     u=\dfrac{\partial u}{\partial \nu}=0 , & \ \ {\rm on}  \ \ \partial \Omega,
     \end{cases}}
\end{equation*}
we have
$$\sum_{i=1}^n(\Gamma_{i+1}-\Gamma_1)^{\frac{1}{2}}\leq
\big(4\Gamma_1^{\frac{1}{2}}+n^2H_0^2\big )^{\frac{1}{2}}\big\{(2n+4)\Gamma_1^{\frac{1}{2}}+n^2H_0^2\big\}^{\frac{1}{2}},
 \eqno{(1.4)}
 $$
where $H_0^2$ is a nonnegative constant which only depends on $M^n$ and $\Omega$. }

\par\bigskip \noindent{\bf Corollary 1}. {\it Under the assumption
of the theorem 1, we have

$$\sum_{i=1}^n\big\{(\Gamma_{i+1}-\Gamma_1)^{\frac{1}{2}}-\Gamma_1^{\frac{1}{2}}\big\}\leq
4\Gamma_1^{\frac{1}{2}}+n^2H_0^2.\eqno{(1.5)}
$$
}

\par\bigskip \noindent{\bf Corollary 2}. {\it When $M^n$ is an $n$-dimensional  complete minimal
submanifold in a Euclidean space, we have
$$
\sum_{i=1}^n(\Gamma_{i+1}-\Gamma_1)^{\frac{1}{2}}\leq
\big\{8(n+2)\Gamma_1\big\}^{\frac{1}{2}}.  \eqno{(1.6)}
 $$
 }
 
\par\bigskip \noindent{\bf Corollary 3}. {\it When $M^n$ is an $n$-dimensional  unit sphere, we have
$$
\sum_{i=1}^n(\Gamma_{i+1}-\Gamma_1)^{\frac{1}{2}}\leq
\big(4\Gamma_1^{\frac{1}{2}}+n^2\big)^{\frac{1}{2}}\big\{(2n+4)\Gamma_1^{\frac{1}{2}}+n^2\big\}^{\frac{1}{2}}.
 \eqno{(1.7)}
 $$
}
\par\bigskip \noindent{\bf Remark 1}. {\it For the unit sphere $S^n(1)$,  by taking  $\Omega=S^n(1)$, 
we know 
$\Gamma_1=0$ and $\Gamma_2=\cdots=\Gamma_{n+1}=n^2$. Hence,
our inequalities become equalities. 
Thus, our results are  sharp.}

\par\bigskip \noindent{\bf Remark 2}. {\it After  the first author and
the third author have proved the above results, 
the second author tells them  that he has also proved the same results . 
Hence, the authors decide   to write this joint paper together.}

\section*{2. Preliminaries}

In order to prove Theorem 1, we need the following Nash's theorem.

\par\bigskip \noindent{\bf Nash's theorem}. {\it Each complete
Riemannian manifold $M^n$ can be isometrically immersed into a
Euclidian space $\mathbb{R}^N$}.

Assume that $M^n$ is an $n$-dimensional
isometrically immersed submanifold in $\mathbb{R}^N$.
Let $\Omega\subset M^n$ be a bounded domain of $M^n$ and $p\in \Omega$. 
Let ($x^1,\cdots,x^n$) be a local coordinate system in a neighborhood $U$ of $p\in M$. Let
${\bf y}$ be the position vector of $p$ in $\mathbb{R}^N$,  which is defined by
$$
{\bf y}=(y^1(x^1,\cdots,x^n),\cdots,y^N(x^1,\cdots,x^n)).
$$
Since $M^n$ is isometrically immersed in $\mathbb{R}^N$, we have
$$g_{ij}=g(\frac{\partial}{\partial x^i}, \frac{\partial}{\partial
x^j})=<\sum\limits_{\alpha=1}^N\frac{\partial y^{\alpha}}{\partial
x^i}\frac{\partial}{\partial y^{\alpha}},
\sum\limits_{\beta=1}^N\frac{\partial y^{\beta}}{\partial
x^i}\frac{\partial}{\partial
y^{\beta}}>=\sum\limits_{\alpha=1}^N\frac{\partial
y^{\alpha}}{\partial x^i}\frac{\partial y^{\alpha}}{\partial
x^j},\eqno{(2.1)}$$ where $g$ denotes the induced metric of $M^n$ from
$\mathbb{R}^N$ and $< , >$ is the standard inner product in
$\mathbb{R}^N$.
The following lemma can be found in \cite{[CC]}.

\par\bigskip \noindent{\bf Lemma}.  {\it For any function $u\in
C^{\infty}(M^n)$, we have
$$\sum\limits_{\alpha=1}^N(g(\nabla y^{\alpha}, \nabla u))^2=|\nabla
u|^2, \eqno{(2.2)}$$
$$\sum\limits_{\alpha=1}^Ng(\nabla y^{\alpha}, \nabla
y^{\alpha})=\sum\limits_{\alpha=1}^N|\nabla
y^{\alpha}|^2=n,\eqno{(2.3)}$$
$$\sum\limits_{\alpha=1}^N(\Delta
y^{\alpha})^2=n^2|H|^2,\eqno{(2.4)}$$
$$\sum\limits_{\alpha=1}^N\Delta y^{\alpha}\nabla
y^{\alpha}=0,\eqno{(2.5)}$$ where $\nabla$ denotes the gradient
operator on $M^n$ and $|H|$ is the mean curvature of $M^n$. }

\section*{3. Proofs of  results}

\vskip 3pt\noindent {\it Proof of  Theorem 1}. Since $M^n$ is a
complete Riemannian manifold, Nash's theorem implies that there
exists an isometric immersion from $M^n$ into a Euclidean space
$\mathbb{R}^N$. Thus, $M^n$ can be considered as an $n$-dimensional
complete isometrically immersed submanifold in $\mathbb{R}^N$.

 Let
$u_i$ be an eigenfunction corresponding to  eigenvalue
$\Gamma_i$ such that $\{u_i\}_{i\in \mathbf{N}}$ becomes an
orthonormal basis of $L^2(\Omega)$, that is,
$$
\int_{\Omega}u_iu_j=\delta_{ij},\ \ \forall i,j\in \mathbf{N}.
$$
We define an $N\times N$-matrix $B$ as follows:
$$B:=(b_{\alpha\beta})$$
where $b_{\alpha\beta}=\int_{\Omega}y^{\alpha}u_1u_{\beta+1}$ and
${\bf y}=(y^{\alpha})$ is the position vector of the immersion in
$\mathbb{R}^N$. Using the orthogonalization of Gram and Schmidt, we
know that there exist an upper triangle matrix $R=(R_{\alpha\beta})$
and an orthogonal matrix $Q=(q_{\alpha\beta})$ such that $R=QB$,
i.e.,
$$
R_{\alpha\beta}=\sum\limits_{\gamma=1}^Nq_{\alpha\gamma}b_{\gamma\beta}=
\int_{\Omega}\sum\limits_{\gamma=1}^Nq_{\alpha\gamma}y^{\gamma}u_1u_{\beta+1}=0,\
\ 1\leq\beta<\alpha\leq N.\eqno{(3.1)}
$$
Defining
$g^{\alpha}=\sum\limits_{\gamma=1}^Nq_{\alpha\gamma}y^{\gamma}$, we
get
$$
\int_{\Omega}g^{\alpha}u_1u_{\beta+1}=\int_{\Omega}\sum\limits_{\gamma=1}^Nq_{\alpha\gamma}
y^{\gamma}u_1u_{\beta+1}=0,\ \ 1\leq\beta<\alpha\leq
N.\eqno{(3.2)}
$$
We put
$$\psi_{\alpha}:=(g^{\alpha}-a^{\alpha})u_1,\ \ \  \ \
a^{\alpha}:=\int_{\Omega}g^{\alpha}u_1^2,\ \ \ \ \ \ 1\leq
\alpha\leq N,\eqno{(3.3)}$$ then it follows that
$$\int_{\Omega}\psi_{\alpha}u_{\beta+1}=0,\ \ \ \
0\leq\beta<\alpha\leq N.\eqno{(3.4)}
$$
Thus, $\psi_{\alpha}$, $1\leq\alpha\leq N$, are trial functions.
From the Rayleigh-Ritz inequality, we have
$$
\Gamma_{\alpha+1}\leq\frac{\int_{\Omega}\psi_{\alpha}\Delta^2\psi_{\alpha}}{\int_{\Omega}\psi_{\alpha}^2},\
\ 1\leq\alpha\leq N.\eqno{(3.5)}
$$
By a direct calculation, we have
$$
\aligned
\int_{\Omega}\psi_{\alpha}\Delta^2\psi_{\alpha}&=\int_{\Omega}\psi_{\alpha}\Delta^2(g^{\alpha}u_1-a^{\alpha}u_1)\\
        &=\int_{\Omega}\psi_{\alpha}\{u_1\Delta^2g^{\alpha}+2\nabla(\Delta
        g^{\alpha})\cdot\nabla u_1+2\Delta g^{\alpha}\Delta u_1\\
        &+2\Delta(\nabla g^{\alpha}\cdot\nabla u_1)+2\nabla
        g^{\alpha}\cdot\nabla(\Delta u_1)+\Gamma_1g^{\alpha}u_1\}.
\endaligned
\eqno{(3.6)}
$$
Then by (3.4) and (3.5), we conclude
$$(\Gamma_{\alpha+1}-\Gamma_1)\|\psi_{\alpha}\|^2\leq\int_{\Omega}r_{\alpha}\psi_{\alpha}:=\omega_{\alpha},
\ \ 1\leq \alpha\leq N, \eqno{(3.7)}$$ where
$$
\aligned
r_{\alpha}&=\Delta^2(g^{\alpha} u_1)-a^{\alpha}\Delta^2u_1\\
           &=u_1\Delta^2g^{\alpha}+2\nabla(\Delta
        g^{\alpha})\cdot\nabla u_1+2\Delta g^{\alpha}\Delta u_1
        +2\Delta(\nabla g^{\alpha}\cdot\nabla u_1)+2\nabla
        g^{\alpha}\cdot\nabla(\Delta u_1).
\endaligned
$$
By making use of Stokes' formula, it is easy to get
$$\int_{\Omega}r_{\alpha}a^{\alpha}u_1=0
$$
and
$$\omega_{\alpha}=\int_{\Omega}r_{\alpha}\psi_{\alpha}=\int_{\Omega}r_{\alpha}(g^{\alpha}u_1-a^{\alpha}u_1)
=\int_{\Omega}r_{\alpha}g^{\alpha}u_1.
$$
We also obtain the following equations from Stokes' theorem
$$2\int_{\Omega}g^{\alpha}u_1\nabla(\Delta g^{\alpha})\cdot\nabla
u_1=\int_{\Omega}\{2u_1\Delta g^{\alpha}\nabla u_1\cdot\nabla
g^{\alpha}+u_1^2(\Delta
g^{\alpha})^2-g^{\alpha}u_1^2\Delta^2g^{\alpha}\},$$
$$2\int_{\Omega}g^{\alpha}u_1\Delta(\nabla g^{\alpha}\cdot\nabla
u_1)=\int_{\Omega}\{2u_1\Delta g^{\alpha}\nabla
g^{\alpha}\cdot\nabla u_1+4(\nabla g^{\alpha}\cdot\nabla
u_1)^2+2g^{\alpha}\Delta u_1\nabla g^{\alpha}\cdot \nabla u_1\},$$
$$2\int_{\Omega}g^{\alpha}u_1\nabla g^{\alpha}\cdot\nabla(\Delta
u_1)=-2\int_{\Omega}\{|\nabla g^{\alpha}|^2u_1\Delta
u_1+g^{\alpha}\Delta u_1\nabla g^{\alpha}\cdot\nabla
u_1+g^{\alpha}\Delta g^{\alpha}u_1\Delta u_1\}.
$$
Consequently, we get
$$
\aligned
\omega_{\alpha}&=\int_{\Omega}\{(\Delta
   g^{\alpha})^2u_1^2+4(\nabla
   g^{\alpha}\cdot\nabla u_1)^2-2|\nabla g^{\alpha}|^2u_1\Delta
   u_1+4u_1\Delta g^{\alpha}\nabla g^{\alpha}\cdot\nabla u_1\}\\
   &=\|u_1\Delta g^{\alpha}+2\nabla g^{\alpha}\cdot\nabla
   u_1\|^2-2\int_{\Omega}|\nabla g^{\alpha}|^2u_1\Delta u_1.
\endaligned
\eqno{(3.8)}
$$
(3.7) and (3.8) imply
$$(\Gamma_{\alpha+1}-\Gamma_1)\|\psi_{\alpha}\|^2\leq \|u_1\Delta g^{\alpha}+2\nabla g^{\alpha}\cdot\nabla
   u_1\|^2-2\int_{\Omega}|\nabla g^{\alpha}|^2u_1\Delta
   u_1, \ \ 1\leq \alpha\leq N.\eqno{(3.9)}$$
On the other hand,
$$
\aligned &\ \ \ \ \int_{\Omega}\psi_{\alpha}(u_1\Delta
            g^{\alpha}+2\nabla u_1\cdot\nabla g^{\alpha})\\
          &=\int_{\Omega}( g^{\alpha}u_1-u_1a^{\alpha})(u_1\Delta
            g^{\alpha}+2\nabla u_1\cdot\nabla g^{\alpha})\\
          &=\int_{\Omega}g^{\alpha}u_1(u_1\Delta
            g^{\alpha}+2\nabla u_1\cdot\nabla g^{\alpha})\\
          &=\int_{\Omega}g^{\alpha}u_1^2\Delta
          g^{\alpha}+\frac{1}{2}\nabla
          u_1^2\cdot\nabla(g^{\alpha})^2.
\endaligned
\eqno{(3.10)}
$$
By using of Stokes' formula, we know
$$\int_{\Omega}g^{\alpha}u_1^2\Delta
  g^{\alpha}=-\int_{\Omega}|u_1\nabla
   g^{\alpha}|^2-\frac{1}{2}\int_{\Omega} \nabla u_1^2\cdot\nabla
   (g^{\alpha})^2.\eqno{(3.11)}
   $$
Substituting (3.11) into (3.10), we infer
$$-\int_{\Omega}\psi_{\alpha}(u_1\Delta
            g^{\alpha}+2\nabla u_1\cdot\nabla g^{\alpha})
=\int_{\Omega}|u_1\nabla
 g^{\alpha}|^2.\eqno{(3.12)}
 $$
From (3.12) and (3.9), we have
$$
\aligned
&(\Gamma_{\alpha+1}-\Gamma_1)^{\frac{1}{2}}\int_{\Omega}|u_1\nabla
     g^{\alpha}|^2\\
     &=-(\Gamma_{\alpha+1}-\Gamma_1)^{\frac{1}{2}}\int_{\Omega}\psi_{\alpha}(u_1\Delta
            g^{\alpha}+2\nabla u_1\cdot\nabla g^{\alpha})\\
            &\leq
            \frac{\delta}{2}(\Gamma_{\alpha+1}-\Gamma_1)\|\psi_{\alpha}\|^2+\frac{1}{2\delta}
            \|u_1\Delta
            g^{\alpha}+2\nabla u_1\cdot\nabla g^{\alpha}\|^2\\
            &\leq (\frac{\delta}{2}+\frac{1}{2\delta})\|u_1\Delta
            g^{\alpha}+2\nabla u_1\cdot\nabla
            g^{\alpha}\|^2 -\delta\int_{\Omega}|\nabla g^{\alpha}|^2u_1\Delta u_1.
\endaligned
\eqno{(3.13)}
$$
According to the lemma in the section 2 and the definition of $g^{\alpha}$, we then have
$$
\aligned &\ \ \ \sum_{\alpha=1}^N\|u_1\Delta g^{\alpha}+2\nabla
  g^{\alpha}\cdot\nabla u_1\|^2\\
&=\sum_{\alpha=1}^N\int_{\Omega}\{u_1^2(\Delta
  g^{\alpha})^2+4(\nabla
  u_1\cdot\nabla g^{\alpha})^2+2(\Delta g^{\alpha}\nabla
  g^{\alpha})\cdot\nabla u_1^2\}\\
&=n^2\int_{\Omega}|H|^2u_1^2+4\int_{\Omega}|\nabla u_1|^2\\
&\leq 4\Gamma_1^{\frac{1}{2}}+n^2\sup\limits_{\Omega}|H|^2.
\endaligned
\eqno{(3.14)}
$$
For any point $p$, by a transformation of coordinates if necessary, we
have, for any $\alpha$,
$$
\aligned
|\nabla g^{\alpha}|^2&=g(\nabla g^{\alpha}, \nabla g^{\alpha})\leq1.
\endaligned
\eqno{(3.15)}
$$

From (3.15), we infer
$$
\aligned
&\sum_{\alpha=1}^N(\Gamma_{\alpha+1}-\Gamma_1)^{\frac{1}{2}}|\nabla
   g^{\alpha}|^2\\
   &\geq\sum_{i=1}^n(\Gamma_{i+1}-\Gamma_1)^{\frac{1}{2}}|\nabla
   g^{i}|^2+(\Gamma_{n+1}-\Gamma_1)^{\frac{1}{2}}\sum\limits_{A=n+1}^N|\nabla g^{A}|^2\\
   &=\sum_{i=1}^n(\Gamma_{i+1}-\Gamma_1)^{\frac{1}{2}}|\nabla
   g^{i}|^2+(\Gamma_{n+1}-\Gamma_1)^{\frac{1}{2}}(n-\sum\limits_{j=1}^n|\nabla
   g^{j}|^2)\\
   &=\sum_{i=1}^n(\Gamma_{i+1}-\Gamma_1)^{\frac{1}{2}}|\nabla
   g^{i}|^2+(\Gamma_{n+1}-\Gamma_1)^{\frac{1}{2}}\sum\limits_{j=1}^n(1-|\nabla
   g^{j}|^2)\\
   &\geq\sum_{i=1}^n(\Gamma_{i+1}-\Gamma_1)^{\frac{1}{2}}|\nabla
   g^{i}|^2+\sum\limits_{j=1}^n(\Gamma_{j+1}-\Gamma_1)^{\frac{1}{2}}(1-|\nabla
   g^{j}|^2)\\
   &=\sum_{j=1}^n(\Gamma_{j+1}-\Gamma_1)^{\frac{1}{2}}.
\endaligned
$$
For (3.13), taking sum on $\alpha$ from $1$ to $N$ and using of (3.14)
and the above inequality, we obtain
$$
\aligned \sum_{i=1}^n(\Gamma_{i+1}-\Gamma_1)^{\frac{1}{2}}
            &\leq
            (\frac{\delta}{2}+\frac{1}{2\delta})
            (4\Gamma_1^{\frac{1}{2}}+n^2\sup\limits_{\Omega}|H|^2)
+n\delta\Gamma_1^{\frac{1}{2}}.
\endaligned
$$
Taking 
$$
\delta=\sqrt{\frac{4\Gamma_1^{\frac{1}{2}}+n^2\sup\limits_{\Omega}|H|^2}
{4\Gamma_1^{\frac{1}{2}}+n^2\sup\limits_{\Omega}|H|^2+2n\Gamma_1^{\frac{1}{2}}}},$$
we obtain
$$
\aligned \sum_{i=1}^n(\Gamma_{i+1}-\Gamma_1)^{\frac{1}{2}}
            &\leq
            (4\Gamma_1^{\frac{1}{2}}+n^2\sup\limits_{\Omega}|H|^2)^{\frac{1}{2}}
            \big\{(2n+4)\Gamma_1^{\frac{1}{2}}+n^2\sup\limits_{\Omega}|H|^2\big\}^{\frac{1}{2}}.
\endaligned
\eqno{(3.16)}
$$
Since the spectrum of the clamped plate problem is an invariant of
isometries, we know that (3.16) holds for any isometric immersion
from $M^n$ into a Euclidean space. 

Now we define $\Phi$ as follows:
$$\Phi:=\{\rm{\psi|\psi\ is\  an\ isometric\ immersion\ from\ {\it{M}} \ into\
a\ Euclidian\ space}\}.
$$
Defining 
$$H_0^2:=\inf\limits_{\psi\in\Phi}\sup\limits_{\Omega}|H|^2,
$$
we obtain

$$ \sum_{i=1}^n(\Gamma_{i+1}-\Gamma_1)^{\frac{1}{2}}\leq
            (4\Gamma_1^{\frac{1}{2}}+n^2H_0^2)^{\frac{1}{2}}
            \big\{(2n+4)\Gamma_1^{\frac{1}{2}}+n^2H_0^2\big\}^{\frac{1}{2}}.
\eqno{(3.17)}
$$
This completes the proof of Theorem 1.

\vskip 5pt\noindent {\it Proof of  Corollary 1}.
Since

$$
            (4\Gamma_1^{\frac{1}{2}}+n^2H_0^2)^{\frac{1}{2}}
            \big\{(2n+4)\Gamma_1^{\frac{1}{2}}+n^2H_0^2\big\}^{\frac{1}{2}}
            \leq
            4\Gamma_1^{\frac{1}{2}}+n^2H_0^2+n\Gamma_1^{\frac{1}{2}}.
$$
from (3.17), we then obtain

$$
\sum_{i=1}^n\big\{(\Gamma_{i+1}-\Gamma_1)^{\frac{1}{2}}-\Gamma_1^{\frac{1}{2}}\big\}\leq
4\Gamma_1^{\frac{1}{2}}+n^2H_0^2. \eqno{(3.18)}
$$
This finishes the proof of Corollary 1.

\vskip 5pt\noindent {\it Proof of  Corollary 2}. For a complete minimal submanifold in
a Euclidean space, we have $|H|=0$. From the proof of the theorem 1,  
the corollary 2 is clear.

\vskip 5pt\noindent {\it Proof of  Corollary 3}. Since  an $n$-dimensional unit sphere 
can be seen as a compact hypersurface  with constant mean curvature $1$ in
the Euclidean space $\mathbb{R}^{n+1}$, we have $|H|=1$. From the proof of the theorem 1,  
the corollary 3 is obvious.

\begin{flushleft}
\medskip\noindent
\begin{tabbing}
XXXXXXXXXXXXXXXXXXXXXXXXXX*\=\kill Qing-Ming Cheng \> Guangyue Huang\\
Department of Mathematics \> Department of Mathematics\\
Faculty of Science and Engineering \> Henan Normal University\\
Saga University \> 453007, Xinxiang \\
840-8502, Saga \> People's Republic of China\\
Japan \>E-mail: hgy@henannu.edu.cn\\
E-mail: cheng@ms.saga-u.ac.jp \\
\end{tabbing}
\end{flushleft}

\begin{flushleft}
\medskip\noindent
\begin{tabbing}
XXXXXXXXXXXXXXXXXXXXXXXXXX*\=\kill Guoxin Wei\\
Department of Mathematics\\
Faculty of Science and Engineering\\
Saga University\\
840-8502, Saga\\
Japan\\
E-mail:wei@ms.saga-u.ac.jp
\end{tabbing}
\end{flushleft}

\end {document}